\documentclass[a4paper,12pt]{article}
\usepackage{amsmath,amsthm}
\numberwithin{equation}{section}

\begin{document}
\centerline{} \centerline{\bf Some Results on $*-$Ricci Flow} \centerline{} \centerline{\bf$~~~~~~$ Nirabhra Basu* and Dipankar Debnath**
}
\centerline{*Department of Mathematics} \centerline{Bhawanipur Education Society College, Kolkata-700020, West Bengal, India}
\centerline{E-mail: nirabhra.basu@thebges.edu.in}
\centerline{**Department of Mathematics}\centerline{Bamanpukur High School(H.S),Nabadwip,India}
\centerline{E-mail: dipankardebnath123@hotmail.com}
\newtheorem{Theorem}{Theorem}[section]
\begin{abstract}\begin{center}
In this paper we have introduced the notion of $*-$ Ricci flow and  shown that $*-$ Ricci soliton which was introduced by Kakimakamis and Panagiotid in 2014, is a self similar soliton of the $*-$ Ricci flow. We have also find the deformation of geometric curvature tensors under $*-$ Ricci flow. In the last two section of the paper, we have found the $\Im$-functional and $\omega-$ functional for $*-$ Ricci flow respectively.
\end{center}\end{abstract}
\section{Introduction}
Ricci flow has been a well known tools of mathematics for last fifteen years after G. Perelman $[1]$ has used it to prove Poincare conjecture in 2003[2]. Sience then the subject has been enriched by many renowned mathematicians.\\\\
A self similar solution to the Ricci flow equation is called Ricci soliton. i.e In differential geometry, a complete Riemannian manifold $(M,g)$ is called a Ricci soliton if and only if there exists a smooth vector field $V$ such that
\begin{equation}
\pounds_V + 2S = \lambda g;
\end{equation}
where $\lambda $ is a constant, $S$ and $\pounds$ denote Ricci curvature tensor$[3]$ and Lie derivative respectively.\\\\
The $*-$ Ricci tensor and $*-$ scalar curvature of $M$ respectively are defined by
$Ric^{*}(X,Y)= \Sigma^{2n+1}_{i=1}R(X,e_i,\phi e_i,\phi Y)$ and $r^{*} = \Sigma^{2n+1}_{i=1} Ric^{*}(e_i,e_i)$ for $X,Y \in TM$
 Tachibana $[7]$  has introduced the notion of $*-$ Ricci tensor on  almost Hermitian manifolds and Hamada $[9]$ use this notation of $*-$ Ricci tensor to almost contact manifold defined by
\begin{equation}
S^{*}(X,Y) = \frac{1}{2}trace(Z\rightarrow R(X,\phi Y)\phi Z)
\end{equation}
for any$~X,Y \in \chi(M);$ where $\chi(M)$ is the Lie algebra of all vector fields on $M$.\\\\
In 2014 Kakimakamis and Panagiotid $[8]$ introduced the concept of $*-$ Ricci solitons within the frame work of real hypersurfaces of complex space form.\\\\
A Riemannian metric $g$ on a manifold $M$ is called a $*-$ Ricci soliton if there exists a constant $\lambda$ and a vector field $V$ such that
\begin{equation}\pounds_Vg+2S^{*}+2\lambda g = 0
\end{equation}
We have defined $*-$ Ricci flow as follows
\begin{equation}
\frac{\partial g}{\partial t}= -2S^*(X,Y)
\end{equation}
In this paper, we have shown that just like Ricci soliton; $*-$Ricci soliton is a self- simmilar soltion of the $*-$ Ricci flow. We have also find the deformation of geometric curvature tensors under $*-$ Ricci flow.\\\\
\textbf{Propostion 1.1: }\\
Defining $\bar{g(t)}=\sigma(t)\phi^{*}_t(g)+\sigma(t)\phi^{*}_t(\frac{\partial g}{\partial t})+\sigma(t)\varphi^{*}_t(\pounds_Xg)$. we have
 \begin{equation}\frac{\partial\hat{g}}{\partial t}= \acute{\sigma}(t)\psi^{*}_t(g)+\sigma(t)+\psi^{*}_t(\frac{\partial g}{\partial t})+\sigma(t)\psi^{*}_t(\pounds_X g)\end{equation}.
 \textbf{Proof:} This follows from the definition of Lie derivative. If we have a metric $g$, a vector field $Y$ and $\lambda\in R $ such that
\begin{equation}-2Ric^{*}(g_0)=\pounds_Yg_{0}-2\lambda g_0 \nonumber\end{equation}
then after setting $g(t)=g_0$ and $\sigma(t) = 1-2\lambda t$ and then integrating the $t-$ dependent vector filed $X(t) = \frac{1}{\sigma(t)}Y. $ To give a family of deffeomorphism $\psi_t $ with $\psi_0$ the identity then $\bar{g}$ defined previously is a Ricci flow with
 \begin{equation}\bar{g} =g_0\frac{\partial\bar{g}}{\partial t} = \sigma'(t)\phi^*_t(g_0)+\sigma(t)\phi^{*}_t(\pounds_Xg_0)\nonumber\end{equation}
 \begin{equation} = \phi^*_t(-2\lambda g_0+ \pounds_Y g_0) = \phi^*_t(-2Ric^*(g_0))=-2Ric^*(\bar{g})\nonumber\end{equation}
 \textbf{Proposition 1.2:}\\
 Under $*-$Ricci  flow
 \begin{equation}g(\frac{\partial}{\partial g}\nabla_XY,Z) = -2(\nabla_XS^*)(Y,Z)+2S^*(Y,\nabla_XZ)+2S^*(\nabla_XY,Z)\nonumber\end{equation}
 \textbf{Proof:}
 $\frac{\partial}{\partial t}\nabla_XY = \pi (X,Y)$
 \begin{equation}g(\frac{\partial}{\partial t}\nabla_XY,Z)=g(\pi(X,Y),Z)\end{equation}
 Again $g(\frac{\partial}{\partial t}\nabla_XY,Z) = \frac{\partial}{\partial t} g(\nabla_XY,Z)- \frac{\partial g}{\partial t}(\nabla_XY,Z)$\\
 \begin{equation}g(\pi(X,Y),Z)= \frac{\partial}{\partial t}g(\nabla_X,Z)+ 2S^*(\nabla_XY,Z)\end{equation}
 Again \begin{equation}Xg(Y,Z) = g(\nabla_XY,Z) + g(Y,\nabla_XZ)\end{equation}
 From $(1.6)$  we have
 \begin{equation}g(\pi(X,Y),Z) = \frac{\partial}{\partial t}[Xg(Y,Z)-g(Y,\nabla_XZ)]+2S^*(\nabla_XY, Z)\nonumber\end{equation}
 \begin{equation}g(\pi(X,Y),Z)= X\frac{\partial g}{\partial t}(Y,Z)-(\frac{\partial g}{\partial t})(Y,\nabla_XZ)+2S^*(\nabla_XY, Z)\nonumber\end{equation}
 or
 \begin{equation}g(\pi(X,Y),Z=-2(\nabla_XS^*)(Y,Z)+2S^*(Y,\nabla_XZ)+2S^*(\nabla_XY, Z)\nonumber\end{equation}
 i.e.
 \begin{equation}g(\frac{\partial}{\partial g}\nabla_XY,Z) = -2(\nabla_XS^*)(Y,Z)+2S^*(Y,\nabla_XZ)+2S^*(\nabla_XY,Z)\end{equation}\section{The $\Im$-functional for the $*$-Ricci flow}
$~~~~~~$Let $M$ be a fixed closed manifold, $g$ is a Riemannian metric and $\textit{f}$ is a function defined on $ M$ to the set of real numbers $\Re$.\\\\
$~~~~~~$ Then the $\Im$-functional on pair $(g,f)$ is defined as\\\\
\begin{equation}\begin{aligned}\Im(g,\textit{f}) = \int( -1 + |\nabla\textit{f} |^{2})e^{-f}dV.\end{aligned}\end{equation}
$~~~~$Now we establish how the $\Im$-functional changes according to time under $*$-Ricci flow.\\

\textbf{Theorem 2.1:}\\
$~~~~~~~$In $*$-Ricci flow, the rate of change of $\Im$-functional with respect of time is given by \\\\
$~~~~\frac{d}{dt}\Im(g,\textit{f})= \int[-2Ric^*(\nabla f,\nabla f)- 2\frac{\partial f}{\partial t}(\Delta f - |\nabla f|^2) + (-1 + |\nabla f|^2)(-\frac{\partial f}{\partial t} + \frac{1}{2}tr\frac{\partial g}{\partial t})]e^{-f} dV$\\\\
where $\Im(g,\textit{f}) = \int( -1 + |\nabla\textit{f} |^{2})e^{-f}dV.$\\\\\\
\textbf{Proof:}
\begin{equation}\begin{aligned}\frac{\partial}{\partial t}|\nabla f|^2 = \frac{\partial}{\partial t}g(\nabla f,\nabla f)
= \frac{\partial g}{\partial t}(\nabla f,\nabla f) + 2g(\nabla \frac{\partial f}{\partial t},\nabla f).\end{aligned}\end{equation}
So using proposition 2.3.12 of $[5]$ we can write\\\\
\begin{equation}\begin{aligned}\frac{d}{dt}\Im( g,f ) = \int[\frac{\partial g}{\partial t}(\nabla f,\nabla f) + 2g(\nabla \frac{\partial f}{\partial t},\nabla f)]e^{-f}dV \\+ \int(-1 + |\nabla f|^2)[-\frac{\partial f}{\partial t} + \frac{1}{2}tr\frac{\partial g}{\partial t}]e^{-f}dV.\end{aligned}\end{equation}
 Using integration by parts of equation$(2.2)$,we get
\begin{equation}\begin{aligned}\int 2g(\nabla \frac{\partial f}{\partial t},\nabla f)e^{-f}dV = -2\int \frac{\partial f}{\partial t}(\Delta f - |\nabla f|^2)e^{-f}dV.\end{aligned}\end{equation}
Now putting (2.4) in (2.3), we get
\begin{equation}\begin{aligned}\frac{d}{dt}\Im( g,f ) = \int[\frac{\partial g}{\partial t}(\nabla f,\nabla f)- 2\frac{\partial f}{\partial t}(\Delta f - |\nabla f|^2)\\ + ( -1 + |\nabla f|^2)(-\frac{\partial f}{\partial t} + \frac{1}{2} tr \frac{\partial g}{\partial t})] e^{-f}dV.\end{aligned}\end{equation}
 Using (1) in (2.5), we get the following result for conformal Ricci flow, as
\begin{equation}\begin{aligned}\frac{d}{dt}\Im(g,f) = \int [-2 Ric^*(\nabla f,\nabla f)  - 2\frac{\partial f}{\partial t}(\Delta f - |\nabla f|^2) \\+ (-1 + |\nabla f|^2)(-\frac{\partial f}{\partial t} + \frac{1}{2}tr\frac{\partial g}{\partial t})]e^{-f}dV .\end{aligned}\end{equation}
Hence the proof.
\section{$\omega$-entropy functional for the $*$- Ricci flow }
$~~~~$ Let $M$ be a closed manifold, $g$ is a Riemannian metric on $M$ and $f$ is a smooth function defined from $M$ to the set of real numbers $\Re$. We define $\omega$-entropy functional as
\begin{equation}\begin{aligned}\omega(g,f,\tau) = \int[\tau(R^* + |\nabla f|^2) + f - n]u dV\end{aligned}\end{equation}
where $\tau > 0 $ is a scale parameter and $ u $ is defined as $ u(t) = e^{-f(t)}~$;$\int_M u dV = 1$.\\\\
We would also like to define heat operator acting on the function $ f:M\times[0,\tau]\longrightarrow\Re$ by
$\diamondsuit:= \frac{\partial}{\partial t}- \Delta $ and also, $\diamondsuit^* := -\frac{\partial}{\partial t} - \Delta +R^*$, conjugate to $\diamondsuit.$\\\\
 We choose $u$, such that $\diamondsuit^* u = 0$.\\\\
 Now we prove the following theorem.\\\\\\
 \textbf{Theorem 3.1:}$~~$ If $g$, $f$, $\tau$ evolve according to
 \begin{equation}\begin{aligned}\frac{\partial g}{\partial t} = -2Ric^*\end{aligned}\end{equation}
\begin{equation}\begin{aligned}\frac{\partial \tau}{\partial t} = -1\end{aligned}\end{equation}
 \begin{equation}\begin{aligned}\frac{\partial f}{\partial t} = -\Delta f + |\nabla f|^2 - R^* + \frac{n}{2\tau}\end{aligned}\end{equation}
 and  the function $v$ defined as $v = [\tau(2\Delta f - |\nabla f|^2 + R^*) + f - n]u,$ then the rate of change of $\omega$-entropy functional for conformal Ricci flow is $~~\frac{d\omega}{dt}= -\int_M \diamondsuit^*v,$ where
\begin{equation}\begin{aligned}\diamondsuit^*v = 2u(\Delta f - |\nabla f|^2 + R^*)-\frac{un}{2\tau} - v - u\tau[ 4<Ric^*,Hess f>\\-2g(\nabla|\nabla f|^2,\nabla f) +4g(\nabla(\Delta f),\nabla f)+ 2|Hess f|^2].\end{aligned}\nonumber\end{equation}
 \textbf{Proof:}
\begin{equation}\begin{aligned}\diamondsuit^*v = \diamondsuit^*(\frac{v}{u}u) = \frac{v}{u}\diamondsuit^*u + u \diamondsuit^*(\frac{v}{u}).\end{aligned}\nonumber\end{equation}
  we have defined previously that$~\diamondsuit^*u = 0,$\\\\
 so$~~ \diamondsuit^*v = u\diamondsuit^*(\frac{v}{u})$\\\\
 $~~~~~~\diamondsuit^*v = u\diamondsuit^*[\tau(2\nabla f - |\nabla f|^2 + R^*) + f - n ].$\\\\
  We shall use the conjugate of heat operator, as defined earlier as$~\diamondsuit^* = -(\frac{\partial}{\partial t} + \Delta - R^*).$\\\\
 Therefore $\diamondsuit^*v = -u(\frac{\partial}{\partial t} + \Delta - R^*)[\tau(2\Delta f - |\nabla f|^2 + R^*) + f -n]$\\\\
 $\Rightarrow u^{-1}\diamondsuit^*v = - (\frac{\partial}{\partial t} + \Delta)[\tau(2\Delta f - |\nabla f|^2 + R^*)] - (\frac{\partial}{\partial t} + \Delta)f - [\tau(2\Delta f - |\nabla f|^2 + R^*) + f - n]$ \\\\
 using equation $(3.3)$, we have
\begin{equation}\begin{aligned}u^{-1}\diamondsuit^*v = (2\Delta f - |\nabla f|^2 + R^*) -\tau (\frac{\partial}{\partial t} + \Delta)(2\Delta f\\- |\nabla f|^2 + R^*) - \frac{\partial f}{\partial t}- \Delta f  - \frac{v}{u}.\end{aligned}\end{equation}
  Now $\frac{\partial}{\partial t}(2\Delta f - |\nabla f|^2 + R^*) = 2 \frac{\partial}{\partial t}(\Delta f) - \frac{\partial}{\partial t}|\nabla f|^2$\\\\
  using proposition $(2.5.6)$ of $[5]$, we have
 \begin{equation}\begin{aligned}\frac{\partial}{\partial t}(2\Delta f - |\nabla f|^2 + R^*) = 2\Delta \frac{\partial f}{\partial t} + 4< Ric^*,Hess f >\\ - \frac{\partial g}{\partial t}(\nabla f,\nabla f) - 2g(\frac{\partial}{\partial t}\nabla f,\nabla f).\end{aligned}\nonumber\end{equation}
 Now using the $*-$Ricci flow equation (1), we have
 \begin{equation}\begin{aligned}\frac{\partial}{\partial t}(2\Delta f - |\nabla f|^2 + R^*) = 2\Delta \frac{\partial f}{\partial t}+~4< Ric^*,Hess f > \\+ 2Ric^*(\nabla f,\nabla f)- 2g(\frac{\partial}{\partial t}\nabla f,\nabla f).\end{aligned}\end{equation}
 Using $(3.4)$ in (3.6), we get
  \begin{equation}\begin{aligned}\frac{\partial}{\partial t}(2\Delta f - |\nabla f|^2 + R^*) = 2\Delta(-\Delta f + |\nabla f|^2 - R^* +\frac{n}{2\tau})\\ + 4<Ric^*,Hess f> + 2Ric^*(\nabla f,\nabla f) - 2g(\frac{\partial}{\partial t}\nabla f,\nabla f).\end{aligned}\end{equation}
 Now let us compute
  \begin{equation}\begin{aligned}\Delta(2\Delta f - |\nabla f|^2 + R^*) = 2\Delta^2f - \Delta |\nabla f|^2.\end{aligned}\end{equation}
 Using $(3.7)$ and $(3.8)$ in $(3.5)$ we obtain after a brief calculation\\\\
$u^{-1}\diamondsuit^*v = (2\Delta f - |\nabla f|^2 + R^*)-\tau[-2\Delta^2f + 2\Delta|\nabla f|^2 + 4 < Ric^*,Hess f > + 2Ric^*(\nabla f,\nabla f)$\\\\
 $~~~~~~~~~~~~$$- 2g(\frac{\partial}{\partial t}\nabla f,\nabla f) + 2\Delta^2f - \Delta|\nabla f|^2 ] - \frac{\partial f}{\partial t}-\Delta f - \frac{v}{u}$\\\\
 $~~~~~~~~~~$$ = \Delta f - |\nabla f|^2 + R^* -\tau[\Delta|\nabla f|^2 + 4< Ric^* ,Hess f> + 2Ric^*(\nabla f,\nabla f)$\\\\
 $~~~~~~~~~ - 2g(\frac{\partial}{\partial t}\nabla f,\nabla f) ] - \frac{\partial f}{\partial t} - \frac{v}{u}$\\\\
 $~~~~~~~~~~~ = \Delta f - |\nabla f|^2 + R^* -\tau[\Delta|\nabla f|^2 + 4< Ric^* ,Hess f> + 2Ric^*(\nabla f,\nabla f)$\\\\
 $~~~~~~~~~~~- 2g(\frac{\partial}{\partial t}\nabla f,\nabla f) ] + \Delta f - |\nabla f|^2 + R^* - \frac{n}{2\tau} - \frac{ v}{u}$\\\\
 $~~~~~~~~~~= 2(\Delta f - |\nabla f|^2 + R^*) - \frac{n}{2\tau} - \frac{v}{u} - \tau[\Delta|\nabla f|^2 + 4<Ric^*,Hess f> + 2 Ric^*(\nabla f,\nabla f)$\\\\
 $~~~~~~~~~~~- 2g(\frac{\partial}{\partial t}\nabla f,\nabla f)]$\\\\
$u^{-1}\diamondsuit^*v = 2(\Delta f - |\nabla f|^2 + R^*)- \frac{n}{2\tau} - [\tau(2\Delta f - |\nabla f|^2 +R^* ) + f - n] - \tau[\Delta|\nabla f|^2$\\\\
$~~~~~~~~~~~~~~~~~~+4<Ric^*,Hess f> + 2 Ric^*(\nabla f,\nabla f)- 2g(\frac{\partial}{\partial t}\nabla f,\nabla f)]$
 \begin{equation}\begin{aligned}u^{-1}\diamondsuit^*v = 2(\Delta f - |\nabla f|^2 + R^*)- \frac{n}{2\tau} - f + n -\tau[2\Delta f - |\nabla f|^2 + R^*\\ + \Delta|\nabla f|^2+ 4<Ric^*,Hess f> + 2Ric^*(\nabla f,\nabla f) - 2g(\nabla \frac{\partial f}{\partial t},\nabla f)]\end{aligned}\end{equation}
using $(3.4)$, we get
 \begin{equation}\begin{aligned}u^{-1}\diamondsuit^*v = 2(\Delta f - |\nabla f|^2 + R^*)- \frac{n}{2\tau} - f + n -\tau[2\Delta f - |\nabla f|^2 \\+ R^* + \Delta|\nabla f|^2+ 4<Ric^*,Hess f> + 2Ric^*(\nabla f,\nabla f)\\- 2g(\nabla(-\Delta f + |\nabla f|^2+\frac{n}{2\tau}-R^*),\nabla f)].\end{aligned}\end{equation}
We can rewrite (3.10) in the following way
 \begin{equation}\begin{aligned}u^{-1}\diamondsuit^*v = 2(\Delta f - |\nabla f|^2 + R^*)-\frac{n}{2\tau} - f + n - \tau[2\Delta f - |\nabla f|^2 +R^* \\+ 4<Ric^*,Hess f> - 2g(\nabla|\nabla f|^2,\nabla f) + 4g(\nabla(\Delta f),\nabla f)]
 \\+ \tau[-\Delta |\nabla f|^2 - 2Ric^*(\nabla f,\nabla f)+2g(\nabla(\Delta f),\nabla f)]\end{aligned}\end{equation}
and using Bochner formula in (3.11) and simplifying, we get\\\\
$u^{-1}\diamondsuit^*v = 2(\Delta f - |\nabla f|^2 + R^*)-\frac{n}{2\tau} - f + n - \tau[2\Delta f - |\nabla f|^2 +R^* + 4<Ric^*,Hess f> $\\\\
$~~~~~~~~~~~~ - 2g(\nabla|\nabla f|^2,\nabla f) + 4g(\nabla(\Delta f),\nabla f)]- 2\tau|Hess f|^2$\\\\
$\Rightarrow u^{-1}\diamondsuit^*v = 2(\Delta f - |\nabla f|^2 + R^*)-\frac{n}{2\tau} - [\tau(2\Delta f - |\nabla f|^2 + R^*)+ f-n]$\\\\$~~~~~~~~~~~~~~~~ - \tau[4<Ric^*,Hess f> -2g(\nabla|\nabla f|^2,\nabla f)+4g(\nabla(\Delta f),\nabla f)]- 2\tau|Hess f|^2$\\\\
 i.e. \begin{equation}\begin{aligned}u^{-1}\diamondsuit^*v = 2(\Delta f - |\nabla f|^2 + R^*)-\frac{n}{2\tau} - \frac{v}{u}- \tau[4<Ric^*,Hess f>\\
 -2g(\nabla|\nabla f|^2,\nabla f)+4g(\nabla(\Delta f),\nabla f)]- 2\tau|Hess f|^2.\end{aligned}\end{equation}
 So finally we have\\\\
 \begin{equation}\begin{aligned}\diamondsuit^*v = 2u(\Delta f - |\nabla f|^2 + R^*)-\frac{un}{2\tau} - v - u\tau[ 4<Ric^*,Hess f>\\
 -2g(\nabla|\nabla f|^2,\nabla f)+4g(\nabla(\Delta f),\nabla f)+ 2|Hess f|^2].\end{aligned}\end{equation}
 Now using remark (8.2.7) of (5), we get \\\\
 $~~~~~~~~~~~~~~~~~~~~~~\frac{d\omega}{dt}= -\int_M \diamondsuit^*v.$\\\\ So the evolution of $\omega$ with respect to time can be found by this integration.
 

\begin{thebibliography}{99}
\bibitem{[1]}Perelman, {\em The entropy formula for the Ricci flow and its geometric applications}, arXiv.org/abs/math/0211159, (2002) 1-39 .
\bibitem{[2]}Perelman, {\em Ricci flow with surgery on three manifolds,} arXiv.org/abs/math/\\0303109, (2002), 1-22.
\bibitem{[3]}R.S. Hamilton, {\em Three Manifold with positive Ricci curvature}, J.Differential
Geom.17(2), (1982), 255-306.
\bibitem{[4]}Bennett Chow, Peng Lu, Lei Ni, {\em Hamilton's Ricci Flow}, American Mathematical Society Science Press, 2006.
\bibitem{[5]}P. Topping, {\em Lecture on The Ricci Flow},Cambridge University Press; 2006.
\bibitem{[6]}A. E. Fischer, {\em An introduction to conformal Ricci flow}, Class. Quantum Grav.21(2004),S171 - S218
\bibitem{[7]}Tachibana S., On almost-analytic vectors in almost Kählerian manifolds, Tôhoku Math. J., 1959, 11, 247–265
\bibitem{[8]}Kaimakamis G., Panagiotidou K., *-Ricci solitons of real hypersurfaces in non-flat complex space forms, J. Geom. Phys., 2014, 86, 408–413
\bibitem{[9]}Hamada T., Real hypersurfaces of complex space forms in terms of Ricci *-tensor, Tokyo J. Math., 2002, 25, 473–483
\end{thebibliography}
\end{document}